\documentclass[twoside,12pt]{article}
\usepackage{graphicx,amsmath,amsfonts,amsthm,amssymb,fancyhdr}

\theoremstyle{definition}
	
\theoremstyle{remark}

\def\beq{\begin{eqnarray}}
\def\eeq{\end{eqnarray}}
\def\bsp{\begin{split}}	
\def\esp{\end{split}}

\newcommand{\be}{\begin{equation}}
\newcommand{\ee}{\end{equation}}

\def \bl {\mbox{\boldmath{$\ell$}}}

\def \hbm #1 {\mbox{\boldmath{$\hat m^{(#1)}$}}}

\def \bm #1 {\mbox{\boldmath{$m^{(#1)}$}}}
\def \BDM {\begin{displaymath}}
\def \EDM {\end{displaymath}}

\begin{document}

\title{$CCNV$ Spacetimes as Potential Supergravity Solutions}
\author{{\large\textbf{A. Coley$^{1}$, D. McNutt$^{2}$ and  N. Pelavas$^{1}$}}
\vspace{0.3cm} \\
$^{1}$Department of Mathematics and Statistics,\\
Dalhousie University,
Halifax, Nova Scotia,\\
Canada B3H 3J5 \\
$^{2}$ Faculty of Science and Technology,\\ 
                         University of Stavanger, 
                         N-4036 Stavanger, Norway \\
\vspace{0.2cm}\\
\vspace{0.3cm} \\
\texttt{aac@mathstat.dal.ca, david.d.mcnutt@uis.no} \\
\texttt{ nicos.pelavas@gmail.com} }
\date{\today}
\maketitle
\pagestyle{fancy}
\fancyhead{}
\fancyhead[CE]{A. Coley, D. McNutt and N. Pelavas}
\fancyhead[LE,RO]{\thepage}
\fancyhead[CO]{$CCNV$ Spacetimes as Potential Supergravity Solutions}
\fancyfoot{}

\begin{abstract}

It is of interest to study supergravity solutions preserving a non-minimal fraction of supersymmetries. A necessary condition for supersymmetry to be preserved is that the spacetime admits a Killing spinor and hence a null or timelike Killing vector field. Any spacetime admitting a covariantly constant null vector field ($CCNV$) belongs to the Kundt class of metrics, and more importantly admits a null Killing vector field. We investigate the existence of additional non-spacelike isometries in the class of higher-dimensional $CCNV$ Kundt metrics in order to produce potential solutions that preserve some supersymmetries. 

\end{abstract}


\begin{section}{Introduction}

Supersymmetric supergravity solutions have been studied in the context of the AdS/CFT conjecture, the microscopic properties of black hole entropy, and the inter-connection of string theory dualities. For example, in five dimensions, solutions preserving various fractions of supersymmetry of $N=2$ gauged supergravity have been studied. The Killing spinor equations imply that supersymmetric solutions preserve $2$, $4,6$ or $8$ of the supersymmetries. The $AdS_{5}$ solution with vanishing gauge field strengths and constant scalars preserves all of the supersymmetries. Half supersymmetric solutions in gauged five dimensional supergravity with vector multiplets possess two Dirac Killing spinors and hence two time-like or null Killing vector fields. These solutions have been fully classified, using the spinorial geometry method, in \cite{gaunt}. Indeed, in  a number of supergravity theories \cite{hommth}, in order to preserve some supersymmetry it is necessary that the spacetime admits a Killing spinor which then yields a null or timelike Killing vector field (isometry) from its Dirac current. Therefore, a necessary (but not sufficient) condition for supersymmetry to be preserved is that the spacetime admits a null or timelike Killing vector field. 

In this short communication we study supergravity solutions preserving a non-minimal fraction of supersymmetries by determining the existence of additional non-spacelike isometries in the class of higher-dimensional Kundt spacetimes admitting a covariantly constant null vector field ($CCNV$) \cite{MCP,Mphd}. $CCNV$ spacetimes belong to the Kundt class because they contain a null Killing vector field which is geodesic, non-expanding, shear-free and non-twisting. The  existence of an additional isometry puts constraints on the metric functions and the vector field components.  Killing vector fields that are null or timelike locally or globally (for all values of the coordinate $v$) are of particular importance. As an illustration we present two explicit examples in this paper. 

A constant scalar invariant ($CSI$) spacetime  is a spacetime such that all of the polynomial scalar invariants  constructed from the Riemann tensor and its covariant derivatives are constant. In three and four dimensions, it has been proven that all $CSI$ spacetimes are either locally homogeneous or belong to the degenerate Kundt class \cite{CSI}; it is conjectured that this is true in higher dimensions as well.  The $VSI$ spacetimes  are $CSI$ spacetimes  for which all of these polynomial scalar invariants vanish. The subset of $CCNV$ spacetimes which are also $VSI$ are of interest. Indeed, it has been shown  previously that the higher-dimensional $VSI$ spacetimes with fluxes and dilaton are solutions of type IIB supergravity \cite{VSISUG}. A subset of Ricci type {\bf N} $VSI$ spacetimes, the higher-dimensional Weyl type\footnote{Weyl and Ricci type are defined in terms of the alignment classification \cite{class}.} {\bf N} pp-wave spacetimes, are known to be solutions in type IIB supergravity with an R-R five-form or with NS-NS form fields \cite{hortseyt, tseytlin}.

In fact, all Ricci type {\bf N} $VSI$ spacetimes are solutions to supergravity   and, moreover, there are $VSI$ spacetime solutions of type IIB supergravity which are of Ricci type {\bf III}, including the string gyratons, assuming appropriate source fields are provided \cite{VSISUG}. It has been argued that the $VSI$ supergravity spacetimes are exact string solutions to all orders in the string tension. Those $VSI$ spacetimes in which  supersymmetry is preserved admit a $CCNV$. Higher-dimensional $VSI$ spacetime solutions to type IIB supergravity preserving some supersymmetry are of Ricci type {\bf N}, Weyl type {\bf III}(a) or {\bf  N} \cite{VSISUG}. 

It is also  known that $AdS_d \times S^{(N-d)}$ spacetimes are supersymmetric $CSI$ solutions of IIB supergravity. There are a number of other $CSI$ spacetimes known to be solutions of supergravity admitting supersymmetries \cite{CSI}, including generalizations of $AdS \times S$ \cite{Gauntlett}, of the chiral null models  \cite{hortseyt}, and the string gyratons \cite{FZ}.  Some explicit examples of $CSI$  $CCNV$  Ricci type {\bf N} supergravity spacetimes have been constructed \cite{VSISUG}. The $CSI$ spacetimes also contain the universal spacetimes \cite{universal}; a spacetime is {\it universal} if every symmetric conserved rank 2 tensor which can be constructed from the metric, the Riemann tensor and its covariant derivatives is proportional to the metric. This ensures that any quantum correction to such a classical solution is proportional to the metric and hence is a solution to all quantum gravity theories defined in terms of a gravitational Lagrangian. Therefore, a subset of the $CSI$ spacetimes will be solutions to all quantum gravity theories defined in this manner. This suggests that the $CSI$ $CCNV$ spacetimes presented in this paper will contain supergravity solutions.

\end{section}

\begin{subsection}{Kundt metrics and $CCNV$ spacetimes} \label{KundtCCNVsect}

A $N$-dimensional spacetime possessing a CCNV, $\ell$, is necessarily of Kundt form. Local coordinates $(u,v,x^e)$ can be chosen, where $\ell = \partial_v$, so that the metric can be written \cite{MCP, Mphd} \begin{equation}
ds^2=2 du [d v+H(u,x^e)d u+ \hat W_{ e}(u,x^f)d x^e]+ {g}_{ef}(u,x^g) dx^e dx^f,  \label{CCNVKundt}
\end{equation} where the metric functions are independent of the light-cone coordinate $v$.

A Kundt metric admitting a  $CCNV$  is $CSI$ if and only if the transverse metric $g_{ef}$ is locally homogeneous \cite{CSI}. Due to the local homogeneity of $g_{ef}$ a coordinate transformation can be performed so that the $m_{ie}$ in equation (\ref{ccnvframe}) below are independent of $u$; this implies that the Riemann tensor is of type {\bf II} or less \cite{class}. If a $CSI$-$CCNV$ metric satisfies $R_{ab}R^{ab}=0$ then the metric is $VSI$, and the Riemann tensor will be of type {\bf III}, {\bf N} or {\bf O} and the  transverse metric is flat (i.e., ${g}_{ef}={\delta}_{ef}$). The constraints on a $CSI$ $CCNV$ spacetime to admit an additional Killing vector field are obtained as subcases of the cases analyzed below where the transverse metric is a locally homogeneous Riemannian manifold.  
\end{subsection}

\begin{section}{$CCNV$ Spacetimes with additional isometries}  \label{GeneralCCNV}
Let us choose the coframe $\{ m^a \}$
\beq m^1 = n = dv + Hdu +  \hat W_e dx^e,~~  \label{ccnvframe} m^2 = \ell,~~ m^i = m^i_{~e} dx^e, \eeq where  $m^i_{~e} m_{if} = g_{ef}$ and $m_{ie}m_j^{~e} = \delta_{ij}$. The frame derivatives are given by:
\beq \ell = D_1 =  \partial_v,~~ n = D_2  =  \partial_u - H\partial_v,~~ m_i = D_i =  m_i^{~e}(\partial_{e} -  \hat W_e \partial_v). \nonumber \eeq 
The Killing vector field can be written as $X = X_1 n + X_2 \ell + X_i m^i$. A coordinate transformation is made to eliminate $ \hat W_3$ in \eqref{CCNVKundt} and we rotate the frame in order to set $X_3 \neq 0$ and $X_m = 0$ \cite{MCP}.  $X$ is now given by
\beq X = X_1 n + X_2 \ell + \chi m^3 . \eeq
Without loss of generality, we assume that the matrix $m_{ie}$ is upper-triangular. As a note, the indices $e,f,g,...$ range from $3$ to $N$ while the indices $m,n,r,p,...$ range from $4$ to $N$.

A subset of the Killing equation  can then be written as:  \beq X_{1,v} = 0,~~  X_{1,u}+X_{2,v} = 0,~~  m_3^{~e}X_{1,e} + X_{3,v} = 0,~~ m_{n}^{\ e}X_{1,e} = 0, \label{letsgo} \eeq which imply 
\beq   & X_1 = F_1(u,x^e),~ X_2 = -D_2(X_1)v+F_2(u,x^e),& \label{kvcomps1} \\~ &X_3 = -D_3(X_1)v +F_3(u,x^e), & \label{kvcomps2}  \eeq and the remaining Killing equations are:
\beq &D_2 X_2  + \displaystyle \sum_i  J_i X_i=0 & \label{killeqn40}\\ & D_i X_2 + D_2 X_i - J_i X_1- \displaystyle \sum_j (A_{ji}+B_{ij})X_j = 0 & \label{killeqn50} \\ & D_j X_i+D_i X_j +2B_{(ij)}X_1 - 2\displaystyle \sum_k \Gamma_{k(ij)} X_k = 0,& \label{killeqn60} \eeq 
\noindent where \beq  & B_{ij} = m_{ie,u}m_{j}^{~e},~ W_i = m_i^{~e} \hat W_e,~~D_{ijk}\equiv  2m_{i e,f}m_{[j}^{~e}m_{k]}^{~f}, & \nonumber \\ &J_i \equiv \Gamma_{2i2} = D_i H - D_2 W_i - B_{ji}W^j,~ A_{ij} \equiv D_{[j} W_{i]}  + D_{k[ij]}W^k,& \\ 
& \Gamma_{ikj} = -\frac12 ( D_{ijk} + D_{jki} - D_{kij}) = {^S}\Gamma_{ikj}. & \nonumber \eeq 
Further information can be found by taking the Killing equations and applying the
commutation relations. This produces two cases; (1) $D_3 X_1 = 0$, or (2) $\Gamma_{3n2} =
\Gamma_{3n3} = \Gamma_{3nm} = 0$.

\begin{subsection}{Case 1: $D_{3}X_{1}=0$}

Using equation \eqref{killeqn40} and the definition of $F_2$ from (\ref{kvcomps1}) and (\ref{kvcomps2}), we have that $X_{1}=c_1u+c_2$. If $c_1 \neq 0$ we may always choose coordinates to set $X_1 = u$, while if $c_1 = 0$ we may choose $c_2 = 1$.   
\vspace{ 3 mm}

\noindent {\bf Subcase 1.1: $F_{3}=0$.} 

\noindent (i) $c_1 \neq 0$,  $X_1 = u$;  $F_2$ must be of the form 
\beq F_2 = \frac{f_2(x^e)}{u} + \frac{g_2(u)}{u}. \label{case11f2} \eeq 
$H$ and $W_m$ are given in terms of these two functions (where $g' \equiv \frac{dg}{du}$) 
\beq H =  \frac{ f_2(x^e) }{ u^2 } - \frac{ g_2'(u) }{ u } +  \frac{ g_2(u) }{ u^2 }, ~~ W_m = \frac{ B_{m}(x^e) }{ u }. \label{C11Wn} \eeq
\noindent (ii) $c_1 = 0$, $X_1 = 1$; ${F_2}_{,u}=0$, and  $H$ and $W_n$ are 
\beq H =  F_2(x^e) + A_0(u, x^r),~~ W_n =  \int D_nA_0 du + C_n(x^e). \label{C11Wn0} \eeq 
\noindent In either case, the only requirement on the transverse metric is that it be independent of $u$. The arbitrary functions in this case are $F_{2}$ and the functions arising from integration.
\vspace{ 3 mm}

\noindent {\bf Subcase 1.2: $F_{3} \neq 0$.}
The transverse metric is now determined by
\beq & m_{33} = -\int \frac{1}{X_1}F_{3,3} du + A_1(x^3, x^r)&, \label{m33,u} \\
& m_{nr,u}= -m_{nr,3} \frac{F_3}{m_{33} X_{1}}, ~~ m_{3r,u}=  - \frac{ F_{3,r}}{X_1} - \frac{m_{3[r,3]}m_3^{~3}F_3}{X_1}. & \label{m3r,u} \eeq
\noindent (i) $c_1 \neq 0$,   $X_1 = u$; $F_i(u,x^e)$ ($i=1,2$) are arbitrary functions, $H$ is given by  
\beq H = - D_2F_2 - \frac{ D_2(F_3^2) }{ 2u } - \frac{ F_3 D_3F_2  }{ u } - \frac{ F_3 D_3(F_3^2)}{ 2u^2 }, \label{C12H} \eeq 
and $W_n$ is determined by 
\begin{equation} D_{2}( u W_n ) + F_3 D_3 W_n + D_n(F_2 - u H ) = 0  . \label{C12Wn} \end{equation}
\noindent (ii) $c_1=0$, $c_2 \neq 0$, $X_1 = 1$;  $F_2$ and $F_3$ satisfy
\beq D_2F_2 + F_3 D_3F_2 + \frac12 D_2(F_3^2) + \frac12 F_3 D_3(F_3^2) = 0. \label{C12f20} \eeq
\noindent $H$ may be written as
\beq H &=& \int m_{33} D_2F_3 dx^3 + F_2 + \frac12 F_3^2 + A_2(u, x^r). \label{C12H0} \eeq
\noindent The only equation for $W_n$ is
\beq & F_3 D_3W_n + D_2 W_n  = D_n(H) . \label{C12Wn0} & \eeq

\noindent (iii) $X_1= 0$:
\vspace{ 3mm} 

We have the following constraints on the functions $m^i_{~e}$:
\begin{equation}
F_{3,3} =0,~~  m_{nr,3} = 0,~~  D_2log(m_{33}) = - \frac{ D_3F_2 }{ F_3 } - D_2log(F_3).   \label{mnr,u0} \end{equation}

\noindent While the metric functions must be of the form: 
\beq & W_n = - \int \frac{ m_{33} D_nF_2}{F_3} dx^3 + E_{n}(u, x^r),~ H - \int \frac{ m_{33} D_2F_2}{F_3} dx^3 + A_3(u, x^r). & \label{C12Wn00} \eeq
There are two further subcases depending upon whether $m_{33,r} = 0$ or not,
whence we may further integrate to determine the transverse metric.
\end{subsection}

\begin{subsection}{Case 2: $\Gamma_{3ia} = 0$} \label{GeneralCCNVc2}
This implies the upper-triangular matrix $m_{ie}$ takes the form: \beq &m_{33} = M_{,3}(u, x^3), m_{3r} = 0, m_{nr} = m_{nr}(u,x^r), & \nonumber \eeq \noindent while the $W_n$ must satisfy $D_3 (W_n) = 0$. The remaining Killing equations  simplify; in particular, $B_{(mn)}X_1 = 0$, leads to two subcases: 1) $X_1 = 0$, or 2) $B_{(mn)} = 0$. 
\end{subsection}
\vspace{ 3 mm}

\noindent {\bf Case 2.1: $X_1 = 0$, $B_{(mn)} \neq 0$.} 
$F_{2,r} = 0$, $F_{3,e} = 0$, with $m_{ie}$,  $H$, $W_n$ given by  equations (\ref{mnr,u0}) and (\ref{C12Wn00}). 
\vspace{ 3 mm}

\noindent {\bf Case 2.2: $B_{(mn)} = 0$, $X_1 \neq 0$.} This case is similar to the subcases dealt with in Case 1.1 (see equations \eqref{case11f2}-\eqref{m33,u}, \eqref{C12Wn0}-\eqref{C12Wn00}). For $n<p$ the vanishing of $B_{(np)}$ implies $m_{n r , u} = 0$, the special form of $m_{ie}$ implies that $m_r^{~~3} = 0$, and the only non-zero component of the tensor $B$ is $B_{33}$. 

If we assume that $F_{1,3} \neq 0$ and $F_1$ is independent of $x^r$: \beq \frac{m_{33,3}}{m_{33}} = \frac{F_{1,33}}{F_{1,3}},~~ \frac{m_{33,u}}{m_{33}} = \frac{F_{1,3u}}{F_{1,3}}. \label{m33comma3} \eeq 
\noindent Thus $m_{33}(u,x^3)$ is entirely defined by $F_1$. We may now solve for $H$ and $W_n$: \beq H = \frac{D_3 D_2 F_1}{D_3(F_1)^2} F_3 - \frac{D_2^2 F_1}{D_3 (F_1)^2} F_1 - \frac{ 2D_{(2 } F_{ 3) } }{ D_3 F_1 },~~ W_n = - \frac{D_n F_3}{D_3F_1}.\label{case22H}   \eeq 
$F_3$ is of the form: \beq F_3 = \int \frac{m_{33} F_1 D_3D_2 F_1 }{D_3F_1} dx^3 + A_6(u, x^r). \label{f3case22} \eeq There are differential equations for $F_2$ in terms  of the arbitrary functions $F_1(u, x^3)$ and $A_6(u, x^r)$. These solutions are summarized in Table 5.2 in \cite{Mphd}. 
\end{section}
\vspace{ 3 mm}

\begin{section}{Killing Lie algebras:}
There are three particular forms for the Killing vector fields in $CCNV$ spacetimes admitting an additional isometry: 
\beq  &(A)& ~~~~ X_A = c n + F_2(u,x^e) \ell + F_3(u,x^e) m^3  \nonumber \\  &(B)& ~~~~ X_B = u n + [F_2(u,x^e)-v] \ell + F_3(u,x^e) m^3 \nonumber \\  &(C)& ~~~~ X_C = F_1(u,x^3) n + [F_2(u,x^e) - D_2F_1 v] \ell + [F_3(u,x^e) - D_3F_1 v] m^3 . \nonumber \eeq 
To determine if these spacetimes admit even more isometries we examine the commutator of $X$ with $\ell$ in each case.  In case (A), $[X_A, \ell ] = 0$ and in case B $[X_B, \ell] = - \ell$, implying there are no additional Killing vector fields.

In the most general case $Y_C \equiv [X_C, \ell]$ can yield a new Killing vector field; $Y_C = D_2F_1 \ell + D_3F_1 m_3$. However,  this will always be spacelike since $(D_3F_1)^2 > 0$.  Note that $[Y_C, \ell ] = 0$, while, in general, $[Y_C, X_C] \neq 0$. 

\begin{subsection}{Globally non-spacelike Killing vector fields}
Let us consider the set of $CCNV$ spacetimes admitting an additional non-spacelike isometry. The equations \eqref{kvcomps1} and \eqref{kvcomps2} imply that the norm of this vector field must satisfy: \beq D_3(X_1)^2  v^2 + 2(D_2(X_1) X_1 - D_3(X_1) F_3) v + F_3^{~2} - 2 X_1 F_2 \leq 0. \nonumber \eeq 

\noindent If the Killing vector field is non-spacelike for all values of $v$,  then $D_3(X_1)$ must vanish  and $X_1$ is constant. Therefore, those subcases with $X_1$ non-constant are excluded.  

We need only consider the Killing vector field of the form $X_A$.
In the timelike case, the subcases with $X_1=0$ are no longer valid as this would imply $F_3^{~2} < 0$. 

In the case that $X_A$ is null. If $c = 0$, $F_3$ must vanish and $F_2$ must be constant, implying that $X$ is a scalar multiple of $\ell$. If $c \neq 0$ we can rescale $n$ so that $2F_2 = F_3^{~2}$, we can then integrate out the various cases: 
\begin{itemize}
\item If $F_3=0$, $F_2$ must vanish as well and the Killing vector field, $X$, is proportional to $n$. The remaining metric functions are now $H = A_0(u,x^r)$ and $W_n = \int D_n(A_0) du + C_n(x^e)$. The transverse metric is unaffected. 
\item If $F_3 \neq 0$, the metric functions must satisfy \beq H = A_2(u,x^r),~D_2(W_n) + D_3(W_n)F_3 = D_n(A_2)\eeq \noindent and $${(log m_{33})_{,u}}= {D_2(log F_3)}.$$ 
\end{itemize} 

In this section we have shown that if we require that the Killing vector field is non-spacelike for all values of $v$, this puts strong conditions on the permitted form of the Killing vector field.  If we only require that the Killing vector field is non-spacelike in a subset of the spacetime there is greater diversity in the choice of Killing vector fields permitted. We present two explicit examples in the following subsections, one which is non-spacelike in a region of the spacetime, and the other which is globally non-spacelike. The $CSI$ $CCNV$ spacetimes admitting globally non-spacelike Killing vector fields, are the above cases where the transverse space is locally homogeneous. 
\end{subsection}

\begin{subsection}{Example 1}
 We first present an explicit example for the case where $X_{1}=u$ and $F_{3}\neq 0$. Assuming that $F_{3}(u,x^i)=\epsilon u m_{33}$ and $\epsilon$ is a nonzero constant, we obtain 
\begin{equation} m_{is,u}+\epsilon m_{is,3}=0 \label{miseqs} \end{equation} 
\noindent and the transverse metric is thus given by 
\begin{equation} m_{is}=m_{is}(x^3-\epsilon u,x^n)\, . \label{missols} \end{equation}
 \noindent We have the  algebraic solution
 \begin{equation} \hat{W}_{3}=-\frac{1}{\epsilon}(H+F_{2,u})-F_{2,3}-\epsilon  m_{33}^{\ \ 2}, \label{w3sol} \end{equation} 
\noindent where $F_{2}(u,x^i)$ is an arbitrary function and $H$  is given by 
\begin{equation} H(u,x^i)=\frac{1}{u}\left[-\int^{u}S(z,x^3-\epsilon u+\epsilon z,x^n) dz + A(x^3-\epsilon u,x^n) \right], \label{hsol} \end{equation} 
\noindent where $A$ is an arbitrary function and $S$ is given by 
\begin{equation} S(u,x^3,x^n)=(uF_{2,u})_{u}+\epsilon u F_{2,3u}+\epsilon^2 u(m_{33}^{\ \ 2})_{u}\, . \end{equation} 
\noindent Furthermore, the solution for $\hat{W}_{n}$, $n=4,\ldots,N$ is 
\begin{equation} \hat{W}_{n}(u,x^i)=\frac{1}{u}\left[-\int^{u}T_{n}(z,x^3-\epsilon u+\epsilon z,x^m) dz + B_{n}(x^3-\epsilon u,x^m) \right] \label{wnsol} \end{equation} 
\noindent where $B_{n}$ are arbitrary functions and $T_{n}$ is given  by 
\begin{equation} T_{n}(u,x^3,x^m)=\left[(uF_{2})_{u}+\epsilon u  F_{2,3}+\epsilon^2 um_{33}^{\ \ 2}\right]_{,n} + \epsilon m_{3n}m_{33}\, . \end{equation}  

In this example, the Killing vector field and its magnitude are given by 
\begin{equation} X=u\mathbf{n}+(-v+F_{2}){\bl}+\epsilon u m_{33}\mathbf{m}^3,~~ X_{a}X^{a}=-2uv+2uF_{2}+(\epsilon u m_{33})^2 \, . \label{Xex1} \end{equation}
\noindent Clearly, the causal character of $X$ will depend on the choice of $F_{2}(u,x^i)$, and for any fixed $(u,x^i)$ $X$ is timelike or null for appropriately chosen values of $v$. Moreover, (\ref{Xex1}) is an example of case (B); therefore the commutator of $X$ and ${\bl}$ gives rise to a constant rescaling of ${\bl}$ and, in general, there are no more Killing vector fields. 

The additional isometry is only timelike or null locally (for a restricted range of coordinate values). However, the solutions can be extended smoothly so that the vector field is timelike or null on a physically interesting part of spacetime. For example, a solution valid on $u>0$, $v>0$ (with $F_{2}<0$), can be smoothly matched across $u=v=0$ to a solution valid on $u<0$, $v<0$ (with $F_{2}>0$), so that the Killing vector field is timelike on the resulting coordinate patch. 

As an illustration, suppose the $m_{3s}$ are separable as follows 
\begin{equation} m_{3s}=(x^3-\epsilon u)^{p_{s}}h_{s}(x^n) \end{equation} 
\noindent and $F_{2}$ has the  form 
\begin{equation} F_{2}=-\frac{\epsilon}{2p_{3}+1}(x^3-\epsilon u)^{2p_{3}+1}h_{3}^{\ 2}+g(u,x^n), \end{equation} 
\noindent where  the $p_{s}$ are constants and $h_{s}$, $g$ arbitrary functions. Thus, from (\ref{hsol}) \begin{equation} H=-\epsilon^2(x^3-\epsilon u)^{2p_{3}-1}[x^3-\epsilon(p_{3}+1)u]h_{3}^{\ 2}-g_{,u} + u^{-1}A(x^3-\epsilon u,x^n),\end{equation}
 \noindent and hence from (\ref{w3sol}) 
\begin{equation}  \hat{W}_{3}=-\epsilon^2p_{3}u(x^3-\epsilon u)^{2p_{3}-1}h_{3}^{\ 2}-(\epsilon  u)^{-1}A(x^3-\epsilon u,x^n). \end{equation} 
\noindent Lastly, equation (\ref{wnsol}) gives \begin{eqnarray} \hat{W}_{n}=\epsilon(x^3-\epsilon u)^{p_{3}}h_{3}\left\{\frac{2(x^3-\epsilon u)^{p_{3}}}{2p_{3}+1}\left[x^3-\epsilon\left(p_{3}+\frac{3} {2}\right)u\right]h_{3,n} \right. & & \nonumber \\ \biggl.\mbox{}-(x^3-\epsilon u)^{p_n}h_{n} \biggr\}- g_{,n}+u^{-1}B_{n}(x^3-\epsilon u,x^m)\,. \end{eqnarray} 
\end{subsection}

\begin{subsection}{Example 2}
A second example corresponding to the distinct subcase where $X_{1}=1$ and assuming  $F_{3}  (u,x^i)=\epsilon m_{33}$ gives the same solutions (\ref{missols}) for the transverse metric (although, in this case, the additional isometry is globally timelike or null). In addition, we have  
\begin{equation} \hat{W}_{3}=\int H_{,3}du + \epsilon^{-1}(F_{2}+f) \label{w3sol2} \end{equation}  \noindent where $H(u,x^i)$, $F_{2} ( x^3-\epsilon u,x^n)$ and $f(x^{i})$ are arbitrary functions. Last, the metric functions $\hat{W}_{n}$ are  
\begin{equation} \hat{W}_{n}(u,x^i)=\int^{u}L_{n}(z,x^3-\epsilon u+\epsilon z,x^m) dz + E_{n} (x^3-\epsilon  u,x^m), \label{wnsol2}\end{equation}
 \noindent with $E_{n}$ arbitrary and $L_{n}$ given by \begin{equation} L_{n} (u,x^3,x^m)=H_{,n}+\epsilon\int H_{,3n}du + f_{,n}\, . \label{Lnsol2} \end{equation} 
\noindent The Killing vector field and its magnitude is 
\begin{eqnarray} X=\mathbf{n}+F_{2}{\bl}+\epsilon m_{33}\mathbf{m}^3,& & X_{a}X^{a}=2F_{2}+(\epsilon m_{33})^2 \, . \label{Xex2} \end{eqnarray}
\noindent Since $F_{2}$ and $m_{33}$ have the same functional dependence there always exists an $F_2$ such that $X$ is everywhere timelike or null. The Killing vector field (\ref{Xex2}) is an example of case (A) and thus $X$ and ${\bl}$ commute and hence no additional isometries arise. For instance, suppose $H=H(x^3-\epsilon u,x^n)$ and $f$ is analytic at $x^3=0$ then (\ref{w3sol2}) and (\ref{wnsol2}) simplify to give 
\begin{eqnarray} \hat{W}_{3} & = & -\epsilon^{-1}(H-F_2 - f), \\ \hat{W}_{n} & = & \epsilon^{-1}\sum^{\infty}_{p=0}\partial_{n}\partial_{3}^{\ p}f(0,x^m)\frac{(x^3)^{p+1}}{(p+1)!} + E_{n} (x^3-\epsilon u,x^m) \, . \end{eqnarray} 

\noindent This explicit solution is an example of a spacetime admitting 2 global null or timelike Killing vector fields, and thus it preserves a non-minimal fraction of supersymmetries. 

\end{subsection}
\end{section}

\begin{section}{Discussion}

The $CCNV$ spacetimes discussed in this paper will preserve a non-minimal fraction of supersymmetries, if they are solutions of some supergravity theory. To show that there are indeed $CCNV$ spacetimes that are solutions to a supergravity theory, we note that there exists $VSI$ and $CSI$ spacetimes which are solutions to supergravity theories  \cite{VSISUG}. The $CSI$ and $VSI$ solutions admit a covariantly constant null vector (i.e., they are $CCNV$-$CSI$ spacetimes) or are constructed from the warped product of a $CCNV$-$VSI$ spacetime and a locally homogeneous Riemannian manifold. 

The construction of $CSI$ solutions for supergravity theories was  motivated by the observation that $AdS_d \times S^{(D-d)}$ is a supersymmetric exact solution of supergravity (for certain values of $(D, d)$ and for particular ratios of the radii of curvature of the two space forms; in particular, $d = 5, D = 10$, $AdS_5  \times S^5$). The more general $D$-dimensional product spacetime $M_d \times K^{(D-d)}$ (in brief $M \times K$) can be considered as a Freund-Rubin background. For example, for $(D, d) = (11, 4),(11, 7)$ and $(5, 5)$ it is sufficient that $M$ and $K$ are Einstein. Since $M \times K$ is a Freund-Rubin background, if $M$ is any Lorentzian Einstein manifold and $K$ is any Riemannian Einstein manifold (with the same ratio of the radii of curvature as in the $AdS \times S$ case), then $M \times K$ will be a solution of some supergravity theory without any consideration of preservation of supersymmetry, where the fluxes are given purely in terms of the volume forms of the relevant factor(s). In general,  $M$ must have negative scalar curvature and $K$ will have positive scalar curvature in order to satisfy the supergravity equations of motion. 

There are many examples of CSI spacetimes in the Freund-Rubin $M \times K$ supergravity set. $K$ could be a homogeneous space or a space of constant curvature. One must ask whether these $CSI$ solutions preserve any supersymmetry. The condition for preservation of supersymmetry demands that $M$ and $K$ admit Killing spinors which imply the existence of Killing vectors. In this paper we have examined the class of spacetimes that will potentially admit more than one Killing spinor without requiring they are solutions to some supergravity theory. Noting that a $CCNV$ spacetime will be $CSI$ if the transverse metric is locally homogeneous we can relate particular instances of the two examples presented in section $3$ to subcases of known $CCNV$-$CSI$ supergravity solutions.

For example, requiring that the transverse space is flat leads to the condition that $m_{33} =1$ and the two $CCNV$ examples in section 3 will contain the subclass of $CCNV$-$VSI$ $\epsilon =0$  spacetimes admitting two timelike or null Killing vectors \cite{VSISUG}-a. The Ricci type ${\bf N}$ $CCNV$-$VSI$ $\epsilon=0$ spacetimes have been shown to be solutions to superstring and heterotic string theory \cite{hortseyt}.
Similarly, imposing the condition that the transverse metric is locally homogeneous will yield $CCNV$-$CSI$  spacetimes admitting two non-spacelike Killing vectors which can be related to subcases of known $CCNV$-$CSI$ solutions of supergravity theories. As a simple example in five dimensions, we may choose the transverse space to be $S^3$ with unit radius, then for an appropriate choice of metric functions, the example in subsection 3.3 corresponds to the metric given by equations (10) and (11) in \cite{VSISUG}-c admitting an additional null or timelike Killing vector field. The metric, together with a constant dilaton and appropriate antisymmetric field is an exact solution of bosonic string theory. 

Motivated by the examples of $CCNV$-$CSI$ spacetimes given in the literature, it is worthwhile to ask if the $CCNV$ spacetimes admitting additional Killing vector fields will contain supergravity solutions beyond the $CSI$ or $VSI$ examples. By removing the condition for the transverse space to be locally homogeneous, such a solution would preserve a non-minimal fraction of supersymmetries while reducing the number of spacelike symmetries. This will be investigated in future work in the context of five-dimensional bosonic string theory.

\end{section}

\begin{section}*{Acknowledgments} 
This work was supported by NSERC of Canada (A.C.) and through the Research Council of Norway, Toppforsk grant no. 250367: Pseudo-Riemannian Geometry and Polynomial Curvature Invariants: Classification, Characterisation and Applications (D.M.).
 
\end{section}


\begin{thebibliography}{99}
\bibitem{gaunt} J.~P.~Gauntlett and J.~B.~Gutowski, Phys.\ Rev.\textit{\ }\textbf{D68} (2003) 105009; J. B.
Gutowski and W. A. Sabra, JHEP \textbf{10}, 039 (2005) \& JHEP \textbf{12 }, 025 (2007) .
\bibitem{hommth} J. M. Figueroa-O'Farrill, P. Meessen and S. Philip, Class. Quant. Grav. {\bf 22}, 207 (2005); E. Hackett-Jones and D. Smith, JHEP {\bf 0411}, 029 (2004).
\bibitem{MCP}  D. McNutt, A. Coley,  and N. Pelavas,   IJGMMP {\bf 6}, 419 (2009).

\bibitem{Mphd} D. McNutt PhD thesis, Dalhousie University (2013). 


\bibitem{CSI}   A. Coley, S. Hervik and N. Pelavas, Class. Quant. Grav. {\bf 23}, 3053 (2006); Class. Quant. Grav.  {\bf 25}  025008 (2007); Class. Quant. Grav. {\bf 26}, 025013 (2009).

 
\bibitem{VSISUG} A. Coley, A. Fuster, S. Hervik and N. Pelavas, JHEP {\bf 32},032 (2007); A. Coley, A. Fuster, S. Hervik and N. Pelavas, Class. Quant. Grav. {\bf 23}, 7431 (2006); A. Coley, A. Fuster and S. Hervik, IJMP A {\bf 24}, 1119 (2009).

\bibitem{class} A. Coley, R. Milson, V. Pravda and A. Pravdov{\'a}, Class. Quant. Grav. {\bf 21} L35 (2004).

\bibitem{hortseyt} G. T. Horowitz  and A. A. Tseytlin, Phys. Rev.  D {\bf 51}, 2896 (1995).
\bibitem{tseytlin} R.~R.~Metsaev and A.~A.~Tseytlin, Phys. Rev. D {\bf 65}, 126004 (2002); J.~G.~Russo and A.~A.~Tseytlin, JHEP {\bf 0209}, 035 (2002); M. Blau {\it et al.}, JHEP {\bf 0201}, 047 (2002); P. Meessen, Phys. Rev. D {\bf 65}, 087501 (2002). 
\bibitem{Gauntlett} J. Gauntlett {\it et al.}, 	Phys. Rev. D{\bf74}, 106007 (2006) .
\bibitem{FZ} V. P. Frolov and A. Zelnikov, Phys. Rev.  D {\bf 72}, 104005 (2005).

\bibitem{universal} A. A. Coley and S. Hervik, ISRN Geometry 2011 (2011).


\end{thebibliography}
\end{document}